\documentclass{article}
\usepackage{amsmath, amssymb, mathrsfs, amsfonts,xfrac}
\usepackage{hyperref}
\usepackage{url}
\usepackage{enumerate,enumitem}
\usepackage{color}
\input xy
\xyoption{all}

\newtheorem{Thm}{Theorem}[subsection]
\newtheorem{Prop}[Thm]{Proposition}

\newtheorem{Lemma}[Thm]{Lemma}
\newtheorem{Cor}[Thm]{Corollary}
\newtheorem{Rmk}[Thm]{Remark}
\newtheorem{Def}[Thm]{Definition}
\newtheorem{Con}[Thm]{Conjecture}
\newtheorem{Question}[Thm]{Question}
\newtheorem{Ex}[Thm]{Example}

\newcommand{\pf}{\noindent\textit{Proof.} }
\newcommand{\qed}{\hfill $\Box$}

\newcommand{\PP}{\mathbb{P}}
\newcommand{\RR}{\mathbb{R}}

\newcommand{\cO}{\mathcal{O}}
\newcommand{\nn}{\mathbf{n}}
\newcommand{\Bir}{\operatorname{Bir}}
\newcommand{\Nef}{\operatorname{Nef}}
\newcommand{\eMov}{\overline{\operatorname{Mov}}^e}
\newcommand{\Mov}{\overline{\operatorname{Mov}}}
\newcommand{\inte}{\operatorname{int}}
\newcommand{\vol}{\operatorname{vol}}
\newcommand{\Eff}{\overline{\operatorname{Eff}}}

\newcommand{\om}{\omega}

\newcommand{\lam}{\lambda}

\newcommand{\inproc}[2]{\langle#1, #2\rangle}
\newcommand{\Abar}{\overline{\mathcal{A}}}
\newcommand{\Cbar}{\overline{\mathcal{C}}}
\newcommand{\Tbar}{\overline{\mathcal{T}}}
\newcommand{\pTbar}{\partial\overline{\mathcal{T}}}
\newcommand{\Fbar}{\overline{\mathcal{F}}}
\newcommand{\CHbar}{\overline{\mathcal{CH}}}
\newcommand{\CH}{\mathcal{CH}}
\newcommand{\ihat}{\widehat{i}}

\newcommand{\khat}{\widehat{k}}

\newcommand{\ww}{\widetilde{w}}
\renewcommand{\l}{\{1,\cdots,l\}}

\begin{document}

\title{Volumes in Calabi-Yau Complete Intersection of Products of Projective Space}
\author{Yi-Heng Tsai}
\date{}
\maketitle

%=========================================================
%=========================================================

\noindent\textbf{Abstract.}
We prove that the birational automorphism group of a general Calabi-yau complete intersection $X$ given by ample divisors in $\PP^{n_1}\times\cdots\times\PP^{n_l}$ is always Lorentzain. Applying the Kawamata-Morrison cone theorem on such $X$, we compute $\vol_X(D+sA)$ for any divisor $D\in \partial\Eff(X)$ and ample divisor $A$ when $s$ is small. We also provide examples of volumes of certain Cartier divisors that involve the digamma function.

\section{Introduction}
It is conjectured in \cite{Mor93} and \cite{kaw97} that the movable effective cone of a Calabi-Yau manifold has a rational fundamental domain under the action of the birational automorphism group.
\begin{Con}(Kawamata-Morrison cone conjecture)
    Let $X$ be a Calabi-Yau manifold. Then there exists a rational polyhedral fundamental domain $\Pi$ for the action of the birational automorphism group $\Bir(X)$ on the movable effective cone $\eMov(X)$ in the sense that \begin{enumerate}
        \item $\eMov(X)=\displaystyle\bigcup_{g \in \Bir(X)} g^*\Pi$,
        \item $\inte(\Pi)\cap \inte(g^*\Pi)=\emptyset$ if $g^* \ne id$.
    \end{enumerate}
\end{Con}
The conjecture has been proven in the case when $X$ is a general Wehler $N$-fold i.e. a general hypersurface of multidegree $(2,\cdots,2)$ in $(\PP^1)^{N+1}$ for $N\ge 3$ (\cite[Theorem 1.3]{CO15}). In \cite{FLT23}, the structure of the boundary of the pesudoeffective cone $\Eff(X)$ has been studied. The divergent–recurrent decomposition of $\partial\Eff(X)$ is proven in \cite[Theorem 2.4.2]{FLT23}, and the following result on the asymptotic behavior of volume function near $\partial\Eff(X)$ is derived:

\begin{Thm}(\cite[Theorem 1.3.7]{FLT23})
    Suppose $X$ is a general Wehler Calabi-Yau $N$-fold. Then, for every pseudoeffective $\RR$-divisor $D$ in $\partial\Eff(X)$ and sufficiently ample divisor $A$, there exists an integer $s(D) \in \{0,\cdots,N-2\}$ and a real number $\delta(D) \in [1,\frac{1}{2}(N-s(D))]$ such that \begin{align*}
        &\liminf_{s \downarrow0}\frac{\log\vol(D+sA)}{\log s}=\delta(D),\\
        &\limsup_{s \downarrow0}\frac{\log\vol(D+sA)}{\log s}=\frac{N-s(D)}{2}.
    \end{align*}
\end{Thm}
The fact that the birational automorphism group of a general Wehler $N$-fold is Lorentzian provides a hyperbolic subspace in $\PP N^1(X)$, and plays an important role when investigating the stucture of $\partial\Eff(X)$ (see \cite[\S 2.3]{FLT23}). Let $\nn=(n_1, \cdots, n_l) \in \mathbb{N}$ such that $|\nn| \ge 4$ and $\nn \ne (2,2)$. In \cite{Yn22}, the author generalized the results in \cite{CO15}, and proved the Kawamata-Morrison cone conjecture when $X$ is a general Calabi-Yau complete intersection in $\PP^\nn:= \PP^{n_1} \times \cdots \times \PP^{n_l}$ given by the intersection of $n$ ample divisors, with $n \le \min{\{n_i\}}$. He also conjectured that the birational automorphism groups for such $X$ are Lorentzian (\cite[Conjecture 4.12]{Yn22}), and gave a nice description of $\partial\Eff(X)$ (\cite[Theorem 1.4]{Yn22}) assuming the conjecture holds. In \S\ref{Ch1}, we will prove this conjecture.
\begin{Thm} \label{Lor}
    Let $X$ be a general Calabi-Yau complete intersection in $\PP^\nn$ defined as the intersection of $n$ ample divisors with $|\nn| \ge 4$ and $\nn \ne (2,2)$. Then $\Bir(X)$ is Lorentzian.
\end{Thm}
As an application of Theorem \ref{Lor}, we will show that the approach in \cite{FLT23} can be applied to such $X$. Following the notation in \cite{FLT23}, we use the notation $A \asymp B$ to indicate that there exists a constant $c>1$ such that $\frac{1}{c}A\le B \le cA$. The asymptotic behavior of volume function near $\partial\Eff(X)$ can be computed as follows:
\begin{Thm} \label{general vol theorem}
    Let $X$ be a general Calabi-Yau complete intersection in $\PP^\nn$ defined by $n$ ample divisors with $|\nn| \ge 4$ and $\nn \ne (2,2)$. Then, for any pseudoeffective divisor $p$ with $\vol_X(p)=0$, any ample $A$, and sufficiently small $s>0$, one of the following holds:
    \begin{enumerate}[label=(\alph*)]
        \item $\vol_X(p+sA)\asymp s^{\frac{kn}{2}}$ for some $k=1 \text{ or } |J|-1$,
        \item $\vol_X(p+sA)\asymp s^{kn}$ for some $k=1, \cdots, |J|-2$,
    \end{enumerate} where $J=\{j\mid n_j=\min\{n_i\}\}$.
    \end{Thm}
    
    It is known that the volume of a Cartier divisor can be written as an integral over the Okounkov body by \cite{LM09}. Recall that a number is called a period or algebraic period if it can be expressed as an integral of an algebraic function over an algebraic domain. Therefore, it is natural to ask\begin{Question}\cite[page 14]{KLM13} \label{Q}
        Is the volume of an integral Cartier divisor on an irreducible projective variety always a period?
    \end{Question}
In \S \ref{Ch4}, considering $\PP^2$-bundles of the general Calabi-Yau complete intersection in $\PP^1\times\PP^1\times \PP^2$, we find the digamma function $\Psi$ naturally arises as the volume of some Cartier divisors.
\begin{Thm} \label{Psi}
    There exists a Cartier divisors $D$ such that $$\vol(D)=\alpha\left(\Psi\left(\frac{7}{4}+\frac{\sqrt{33}}{12}\right)-\Psi\left(\frac{7}{4}-\frac{\sqrt{33}}{12}\right)+\Psi\left(\frac{5}{4}+\frac{\sqrt{33}}{12}\right)-\Psi\left(\frac{5}{4}-\frac{\sqrt{33}}{12}\right)\right)+\beta$$ where $\alpha$ and $\beta$ are non-zero periods.
\end{Thm}
Since all periods form a ring, if Question \ref{Q} is true, then this would imply the number $\Psi\left(\frac{7}{4}+\frac{\sqrt{33}}{12}\right)-\Psi\left(\frac{7}{4}-\frac{\sqrt{33}}{12}\right)+\Psi\left(\frac{5}{4}+\frac{\sqrt{33}}{12}\right)-\Psi\left(\frac{5}{4}+\frac{\sqrt{33}}{12}\right)$ is a period. However, at most irrational points, the digamma function $\Psi$ has neither been expressed in closed form nor proven to be a period. Therefore, the theorem above is a potential counterexample to Question \ref{Q}.
\\

\noindent \textbf{Acknowledgements} The author is deeply grateful to Cristopher Hacon for his valuable suggestions that greatly improved the quality of this work. The author also thanks José Yáñez for useful discussions related to Kawamata-Morrison cone conjecture. The author was partially supported by  NSF research grant DMS-2301374 and by a grant from the Simons Foundation SFI-MPS-MOV-00006719-07. The author also acknowledges the partial inspiration derived from a visit to the National Center for Theoretical Sciences (NCTS), Taiwan.

\section{Calabi-Yau complete intersection of Products of Projective space}
\label{Ch1}
Let $\nn=(n_1, \cdots, n_l) \in \mathbb{N}^l$ such that $|\nn| \ge 4$ and $\nn \ne (2,2)$. Let $X$ be a general complete intersection Calabi-Yau subvariety in $\PP^\nn:= \PP^{n_1} \times \cdots \times \PP^{n_l}$ given by the intersection of $n$ ample divisors, with $n = \min{\{n_i\}}$, then the Kawamata-Morrison conjecture holds for such $X$ \cite[Theorem 1.3]{Yn22}. In the following section, we'll describe the action of $\Bir(X)$ on the effective movable cone $\eMov(X)$.

\begin{Rmk}
    See \cite[Remark 3.10 and 4.8]{Yn22} for the case when $n < \min{\{n_i\}}$.
\end{Rmk}

\subsection{Kawamata-Morrison cone conjecture}

Let $\pi_i:X\longrightarrow \PP^{n_1}\times \cdots \times \widehat{\PP^{n_i}}\times \cdots \times \PP^{n_l}$ be the natural projections. Then $\pi_i$ defines a birational involution $\iota_i$ for $i\in J:=\{j|n_j=n\}$ \cite[Proposition 2.6]{Yn22}. Furthermore, \cite[Proposition 3.6]{Yn22} gives a matrix description of the action of $\iota^*_j$ on $N^1(X)$ for $j\in J$:

\begin{Prop} \label{iota}
    Let $H_i$ be the divisor class of the pullback of $\cO_{\PP^{n_i}}(1)$ via the natural projection $\PP^{\nn} \longrightarrow \PP^{n_i}$, and $h_i$ be the pullback of $H_i$ via the inclusion $X \longrightarrow \PP^{\nn}$. Then $\iota^*_j (h_i)=h_i$ if  $i\ne j$ and $$
        \iota^*_j (h_j)=-h_j+ \sum_{i\ne j} b_{ij} h_i$$ where $b_{ij}$ is defined to be the coefficient of $H_iH_j^{n-1}$ in $X$. In other words, the matrix of $\iota^*_j$ with respect to the basis $\{h_i|i =1, \cdots ,l\}$ is $$\begin{bmatrix}
        1  &\cdots & 0 & b_{1j} &0 & \cdots &0\\
        \vdots  &\ddots & \vdots & \vdots &\vdots  &\ddots & \vdots\\
        0 &\cdots &1 &b_{j-1j} &0 & \cdots &0\\
        0 &\cdots &0 &-1 &0 & \cdots &0\\
        0 &\cdots &0 &b_{j+1j} &1  &\cdots & 0 \\
        \vdots  &\ddots & \vdots & \vdots &\vdots  &\ddots & \vdots\\
        0 & \cdots &0 &b_{lj}&0 &\cdots &1
    \end{bmatrix}$$
\end{Prop}

Let $W_J$ be the subgroup of $\Bir(X)$ generated by $\iota^*_j$ for all $j \in J$, then the $\Bir(X)$-orbit (of a divisor $D$) is equivalent to the $W_J$-orbit (of $D$) \cite[Proposition 3.3 and 3.8]{Yn22}, and we have \cite[Theorem 1.3]{Yn22}:
\begin{Thm}\label{KMcone}
    $\Nef(X)$ is a fundamental domain for the action $\Bir(X)$, or equivalently $W_J=\langle\iota^*_j|j\in J, \iota^{*2}_j=id\rangle$, on $\eMov(X)$ in the sense that \begin{enumerate}
        \item $\eMov(X)=\displaystyle\bigcup_{g \in \Bir(X)} g^*\Nef(X)$,
        \item $\inte(\Nef(X))\cap \inte(g^*\Nef(X))=\emptyset$ if $g^* \ne id$.
    \end{enumerate}
\end{Thm}

\subsection{The boundary of $\Mov(X)$}
In \cite{Yn22}, the author gave a nice description of the boundary of $\Mov(X)$ under the assumption that the birational automorphism group $\Bir(X)$ is Lorentzian \cite[Theorem 1.4]{Yn22}, and conjuctured $\Bir(X)$ is always Lorentzian \cite[Conjecture 4.12]{Yn22}. In the following, we will introduce some necessary notations and definitions, and then prove $\Bir(X)$ is always Lorentzian for a general Calabi-Yau complete intersection $X$ in $\PP^\nn$.
\begin{Def} \label{B} (cf. \cite[Section 2.1 and Proposition 3.5]{Yn22})
    The birational automorphism group $\Bir(X)$, or equivalently $W_J$, is said to be Lorentzian if the corresponding $|J|\times|J|$ symmetric matrix $B_J$ defined below has signature $(|J|-1,1)$. $$\begin{cases}
       (B_J)_{jj}= 1 \quad \text{ for all } j=1, \cdots, |J|,\\
       (B_J)_{jk}=-\displaystyle\frac{b_{jk}}{2}=-\frac{b_{kj}}{2} \quad\text{ if } j \ne k \le |J|.
    \end{cases}$$
\end{Def}

\begin{Rmk} \label{bij remark}
    By the proof of \cite[Proposition 3.5]{Yn22}, we know that $b_{jk}$ must be either $2n$ or $2n+1$ fro all $j \ne k \le |J|$.
\end{Rmk}

We now provide a proof of Theorem \ref{Lor}. 
\begin{Prop} \label{Lorentzian}
    Let $l>1$ and $0=r_0<r_1\le\cdots\le r_k= l$ be a partition of $l$. Let $A_\alpha$ be the $\alpha\times \alpha$ matrix given by \begin{align*}
    (A_\alpha)_{ij}= \begin{cases}
        1 &\text{if } i=j,\\
        -n  &\text{if } i\ne j
    \end{cases}
\end{align*} where $n \ge 1$. Let $B$ be the matrix formed by the blocks $A_{\alpha_i}$ with $\alpha_i:=r_i-r_{i-1}$ on the diagonal, and outside of it the entries of the matrix only have the value $-n-\frac{1}{2}$. Then the signature of $B$ is $(l-1,1)$.
\\\\
$B=
\left [
    \begin{array}{ccccc|c|ccc}
     1 & -n & -n & \cdots & -n &&&&\\
    -n & 1 & -n & \cdots & -n &&&&\\
    \vdots & & \ddots & & \vdots &-n-\frac{1}{2}&&-n-\frac{1}{2}&\\
    -n & \cdots & -n & 1 & -n &&&&\\
    -n &-n & \cdots & -n & 1 &&&&\\
    \hline
    &&-n-\frac{1}{2}&&& \ddots &&-n-\frac{1}{2}&\\
    \hline 
    &&&&&&1 & \cdots & -n \\
    &&-n-\frac{1}{2}&&&-n-\frac{1}{2}&\vdots & \ddots &  \vdots\\
    &&&&&&-n & \cdots & 1\\
\end{array}\right ]
\begin{array}{c}
         \\\\\alpha_1\\\\\\\vdots\\\\\alpha_k\\\\ 
    \end{array}\\
    \begin{array}{ccc}
        \hspace{3.25cm}\alpha_1&\hspace{2.25cm}\cdots&\hspace{1.75cm}\alpha_k
    \end{array}$
\end{Prop}
\begin{Rmk}
    According to \cite[Proposition 4.13]{Yn22}, the preceding proposition is equivalent to \cite[Conjecture 4.12]{Yn22}.
\end{Rmk}

\textit{Proof.} Let $L_{ij}(m)$ be the elementary matrix corresponding to adding $m$ times row $j$ to row $i$, namely, the left multiplication by $L_{ij}(m)$ adds $m$ times row $j$ to row $i$. Then the right multiplication by $L_{ij}(m)^{-1}=L_{ij}(-m)$ adds $-m$ times column $i$ to column $j$. Define $B':=\left(\prod_{i=1}^k \prod_{j=r_{i-1}+2}^{r_i} L_{j,r_{i-1}+1}(-1)\right)\cdot B\cdot \left(\prod_{i=1}^k \prod_{j=r_{i-1}+2}^{r_i} L_{j,r_{i-1}+1}(-1)\right)^{-1}$. Then, the entries of $B'$ are given by
$$B'_{pq} = \begin{cases}
    1-(r_i-1)n &\text{if } p=q=r_{i-1}+1,\\
    1+n &\text{if } p=q\ne r_{i-1}+1,\\
    -(n+1/2)r_i &\text{if } p=r_{j-1}+1,q=r_{i-1}+1,p\ne q,\\
    -n &\text{if } p=r_{j-1}+1,q=r_{j-1}+2,\cdots, r_j,\\
    -n-1/2 &\text{if } p=r_{j-1}+1,q\ne r_{j-1}+1,\cdots, r_j,q\ne r_{i-1}+1\\
    0 &\text{otherwise.}
\end{cases}$$

Let $B''$ be the matrix obtained from $B'$ by interchanging row $i$, row $r_{i-1}+1$, and interchanging column $i$, column $r_{i-1}+1$ for $i=2,\cdots, k$. Note that interchanging rows or columns does not change the absolute value of the determinant. Consequently, combining with the fact that $B$ is similar to $B'$, we conclude that the eigenvalues of $B$ are $n+1$ with multiplicity $l-k$, and the remaining eigenvalues are those of $A$ where is the submatrix obtained from $B''$ by deleting rows $k+1,\cdots,l$ and columns $k+1,\cdots,l$.

We are now reduced to showing that $A$ has signature $(k-1,1)$. By direct computation, we find
\begin{align*}
    &L_{k1}(-1)\cdots L_{21}(-1)\cdot A\cdot L_{21}(1)\cdots L_{k1}(1)\\
&= \begin{bmatrix}
    1-(\alpha_1-1)n-\sum_{i=2}^k (n+\frac{1}{2})\alpha_i & -(n+\frac{1}{2})\alpha_2 & \cdots & \cdots& -(n+\frac{1}{2})\alpha_k\\
    \frac{1}{2}(\alpha_2-\alpha_1)& 1+n+\frac{\alpha_2}{2} &0 &\cdots &0\\
    \vdots & 0 & \ddots && \vdots\\
    \vdots & \vdots && \ddots & 0 \\
    \frac{1}{2}(\alpha_k-\alpha_1) & 0 &\cdots & 0& 1+n+\frac{\alpha_k}{2}
\end{bmatrix}
\end{align*}
Using the Laplace expansion along the first row, we compute the characteristic polynomial \begin{align*}
    \det{(A-tI)} =& \left(1-(\alpha_1-1)n-\sum_{i=2}^k (n+\frac{1}{2})\alpha_i - t\right)\prod_{i=2}^k (1+n+\frac{\alpha_i}{2}-t)\\&+\sum_{i=2}^k(n+\frac{1}{2})\alpha_i\cdot \frac{1}{2}(\alpha_i-\alpha_1)\prod_{j \ne 1,i}(1+n+\frac{\alpha_j}{2}-t).\\
\end{align*}
Substituting $t=0$ to compute the determinant of $A$, we obtain \begin{align*}
     \det{A} =&\left(1-(\alpha_1-1)n-\sum_{i=2}^k (n+\frac{1}{2})\alpha_i\right)\prod_{i=2}^k (1+n+\frac{\alpha_i}{2})\\&+\sum_{i=2}^k(n+\frac{1}{2})\alpha_i\cdot \frac{1}{2}(\alpha_i-\alpha_1)\prod_{j \ne 1,i}(1+n+\frac{\alpha_j}{2})\\
     =& (1-(\alpha_1-1)n)\prod_{i=2}^k (1+n+\frac{\alpha_i}{2}) - \sum_{i=2}^k (n+\frac{1}{2})\alpha_i\left(\prod_{j\ne i}1+n+\frac{\alpha_j}{2}\right).
\end{align*}
Now, assume $\alpha_1 \le \cdots \le \alpha_k$, then $$\det{(A-tI)}=\begin{cases}
    \det{A<0} \quad\text{ if } t=0,\\
    (n+\frac{1}{2})\alpha_i\cdot \frac{1}{2}(\alpha_i-\alpha_1)\prod_{j \ne 1,i}(\frac{\alpha_j-\alpha_i}{2})\quad\text{ if } t=1+n+\frac{\alpha_i}{2} \text{ for } i=2,\cdots,k
\end{cases}$$
Note that as $t\to-\infty$, $\det{(A-tI)}\to+\infty$. By the intermediate value theorem, the matrix $A$ has at least one negative eigenvalue, and at least $k-1$ positive eigenvalues. Since the total number of eigenvalues of $A$ is $k$, $A$ has exactly one negative eigenvalue and $k-1$ positive eigenvalues, giving $A$ the signature $(k-1,1)$. This completes the proof.\qed

\begin{Cor} \label{invBdiag}
    All diagonal entries of $(B_J)^{-1}$ are positive.
\end{Cor}
\pf This is a direct consequence of Cramer's rule. \qed
\\

We can then describe the boundary of $\Mov(X)$.
\begin{Thm} \label{boundary}
    The boundary of $\Mov(X)$ is the closure of the union of:\begin{enumerate}
        \item The $\Bir(X)$-orbit of the codimension one faces $\{\sum_{k\ne j}a_kh_k|a_k\ge 0\}$ for $j\in J$
        \item The $\Bir(X)$-orbit of the cones $\{a_\lambda v_\lambda+\sum_{k \ne i,j}a_kh_k|a_k\ge 0, a_\lambda \ge 0\}$ for $i,j \in J$ where $v_\lambda$ is an eigenvector associated to the unique eigenvalue $\lambda>1$ of $(\iota_i\iota_j)^*$ if $n\ge 2$; or $v_\lambda=0$ if $n=1$. 
    \end{enumerate}
\end{Thm}

\pf The theorem follows directly from Proposition \ref{Lorentzian} and \cite[Theorem 1.4]{Yn22}. \qed

\begin{Rmk}
    See \cite[\S 5]{Yn22} for applications of Proposition \ref{Lorentzian} to the numerical dimension $v_{\vol}^\RR$.
\end{Rmk}

\section{Volume functions}
Thanks to Theorem \ref{KMcone}, we have a sufficient number of birational automorphisms to control the pseudoeffective cone $\Eff(X)$, with the nef cone $\Nef(X)$ serving as a rational fundamental domain. Consequently, it is natural to expect a clear understanding of the behavior of volume functions on the boundary of $\Eff(X)$ \cite[\S1.6]{FLT23}. 

In \cite{FLT23}, assuming $X$ is a Wehler manifold, i.e. a general hypersurface in $\PP^1\times \cdots \PP^1$ of degree $(2, \cdots ,2)$, they investigated the function $s \mapsto \vol_X(D+sA)$ for a pseudoeffective divisor $D$ with $\vol_X(D)=0$ and any sufficiently ample $A$. We will show the approach in \cite{FLT23} can be applied to a general complete intersection Calabi-Yau subvariety in $\PP^{n} \times \cdots \times \PP^{n}$, the product of $l$ copies of $\PP^n$. For the rest of this paper, we assume that $nl \ge 4$, and $l \ge 3$ if $n=2$.

\subsection{Preliminaries} \label{preliminaries}
We will extend the concept introduced in \cite[\S 2]{FLT23} for the case $n=1$ to the more general case where $n \ge 1$. Throughout this section, let $\iota_i$ denote the birational involution induced by the natural projection $X \to \PP^{n} \times \cdots \times \widehat{\PP^{n}} \times \cdots \times \PP^{n}$ (omitting the $i$th factor $\PP^n$), and let $h_i$ represent the pullback of $\cO_{\PP^{n}}(1)$ via the natural projection $X \to \PP^{n}$ corresponding to the $i$th $\PP^n$. Also, we will abbreviate the matrix $B_J$ as $B$, and $W_J$ as $W$ for simplicity.
\\

\noindent\textbf{The bilinear form.} We define the bilinear form on $N^1(X) \times N^1(X)$ corresponding to $-B^{-1}$ (cf. \cite[2.2.2]{CO15} and \cite[2.2.3]{FLT23}): \begin{equation}\label{hihj}
    \inproc{h_i}{h_j}:=-B_{ij}^{-1}
\end{equation} 
Then the dual basis $\alpha_j$ with respect to $\inproc{\quad}{\quad}$, i.e., $\inproc{h_i}{\alpha_j}=\delta_{ij}$, is given by \begin{equation} \label{alpha}
    \alpha_j= -\sum_k B_{kj}h_k=-h_j+\sum_{i\ne j} \frac{b_{ij}}{2}h_i,
\end{equation} and $\inproc{\alpha_i}{\alpha_j}=-B_{ij}$ by direct computation. Note that Proposition \ref{iota} implies $\iota^*_j$ is the orthogonal reflection in $\alpha_j$: \begin{equation} \label{iota_alpha}
    v \mapsto v-2\frac{\inproc{v}{\alpha_j}}{\inproc{\alpha_j}{\alpha_j}}\alpha_j = v+2\inproc{v}{\alpha_j}\alpha_j.
\end{equation}

\begin{Lemma} \label{equi}
    The action of $W$ preserves the bilinear form $\inproc{\quad}{\quad}$, namely, $\inproc{wv}{wu}=\inproc{v}{u}$ for all $w\in W$ and $u,v\in N^1(X)$.
\end{Lemma}
\pf We may assume $w=\iota^*_j$, $v=h_i$, and $u=h_k$, with $i,j,k \in \l$. Using (\ref{iota_alpha}), we compute $\inproc{\iota^*_jh_i}{\iota^*_jh_k} = \inproc{h_i+2\delta_{ij}\alpha_j}{h_k+2\delta_{jk}\alpha_j} = \inproc{h_i}{h_k}+4\delta_{ij}\delta_{jk}(1-B_{jj})$. The lemma follows immediately from the fact that $B_{jj}=1$.\qed
\\

\noindent\textbf{Cones.} Following the notations in \cite[\S 2]{FLT23}, we let $\Abar$ be the closed cone generated by $\{h_i\}_{i=1,\cdots,l}$, and let $\Cbar=\{v\in N^1(X)| \inproc{v}{a}\ge 0 \text{ for all } a\in\Abar\}$ be its dual cone. Then $\Abar$ if the nef cone of $X$ by \cite[Proposition 2.5]{Yn22}. Consider the $W$-orbits of $\Abar$, we define the {\it Tits cone} $\mathcal{T}:= W\cdot\Abar$. By Theorem \ref{KMcone}, we know that $\mathcal{T}$ is the movable effective cone of $X$. Furthermore, as a direct consequence of \cite[Lemma 4.5]{Yn22}, its closure $\Tbar$ coincides with the pseudoeffective cone of $X$. Let $\CHbar$ be the dual cone of $\Tbar$, then the fundamental domain for the action $W$ on the interior 
$\CH$ of $\CHbar$ is contained in $\Fbar:=\Abar\cap\Cbar$ by \cite[Theorem 5.2]{DLK}.
\\

\noindent\textbf{Limit roots and limit sets.} We recall the notions introduced in \cite[\S 2]{CL17} and \cite[\S2.3]{FLT23} for subsequent use. Let $\Tbar_+$ be the set of equivalence classes of $x\in \Tbar$ under the equivalence relation given by $x\sim \lam x$ for positive $\lam$, and let $\mathcal{T}_+$ denote the image of $\mathcal{T}$ under $\Tbar \longrightarrow \Tbar_+$. For any subset $S \subset \{1,\cdots,l\}$, a limit root $\hat{x}$ of $W_S:=\langle\iota^*_s\mid s\in S\rangle$ is defined as a limit point of $W_S$ acting on $\mathcal{T}_+$. In other words, there exists a base point $\hat{x}_0\in \mathcal{T}_+$ and an injective sequence $w_k \cdot \hat{x}_0$ in the orbit $W_S \cdot  \hat{x}_0$ such that $\displaystyle\hat{x}=\lim_{k\to \infty} w_k \cdot \hat{x}_0$. We define the limit set $\Lambda_S \subset \Tbar_+$ of $W_S$ as the set of all such limits of $W_S$. By \cite[Theorem 12.1.3]{Ra19}, $\Lambda_S$ is contained in $\CHbar$.

\begin{Thm}\cite[Corollary 2.7]{CL17}
    Limit roots are independent of the choice of base point.
\end{Thm}

\begin{Thm}\cite[Theorem 2.7(ii)]{HLR14}
    For any $\hat{x}\in \Lambda_{\l}$, $\inproc{\hat{x}}{\hat{x}}=0$.
\end{Thm}

\noindent\textbf{Codimension 1 decomposition.} Let $S\subset \l$, and $X'$ denote the restriction of $X$ obtained through the natural embedding of $\displaystyle\left(\prod_{i\in S}\PP^n\right) \times \left(\prod_{i\in S^c}\{\text{point}\}\right)$ into $\displaystyle\prod_{i\in \l}\PP^n$. Then $N^1(X')$ is generated by $h_i'$ for $i\in S$, where $h_i'$ denotes the pullback of $\cO_{\PP^{n}}(1)$ via the natural projection $X' \to \PP^{n}$ corresponding to the $i$th $\PP^n$. We then view $N^1(X')$ as a subspace of $N^1(X)$ by identifying $h_i'$ with $h_i$. Similarly, we can define the subgroup $W'=\langle \iota_j^{'*}\rangle$ of $\Bir(X')$, as well as the cones $\Abar'$, $\Tbar'$, and other related objects, analogous to their corresponding definitions for $X$. To prove Theorem \ref{decomp} by induction on $\dim (N^1(X))$, it is essential to establish the relationship between $\Tbar'$ and $\Tbar$ when $|S^c|=1$.
\begin{Prop} (cf. \cite[Proposition 2.4.7]{FLT23})
    Assume $l\ge 3$ and $S^c=\{l\}$. Define the orthogonal projection $$\phi: N^1(X) \longrightarrow N^1(X) \quad\text{ by } v \mapsto v-\frac{\inproc{v}{h_l}}{\inproc{h_l}{h_l}}h_l.$$
    \begin{enumerate}
        \item The map $\phi$ respects the action of $W'$ and $\langle\iota^*_j\rangle_{j=1,\cdots,l-1}$ on $N^1(X')$, meaning that for any $w \in W'$ and $v \in N^1(X')$, $\phi(w \cdot v) = \phi_W(w) \cdot \phi(v)$, where $\phi_W:W'\longrightarrow \langle\iota^*_j\rangle_{j=1,\cdots,l-1} \subset W$ is defined by $\iota_j^{'*}\mapsto\iota^*_j$.
        \item $\Tbar\cap h_l^\perp=\phi(\Tbar')$ where $h_l^\perp$ denotes the set $\{v\in N^1(X)\mid \inproc{v}{h_l}=0\}$.
        \item Let $h_l^{\le 0}=\{v\in N^1(X)\mid \inproc{v}{h_l}\le 0\}$. Then $\Tbar\cap h_l^{\le 0}$ is the convex hull of $\RR_{\ge 0}\cdot h_l \cup (\Tbar\cap h_l^\perp)$.
    \end{enumerate}
\end{Prop}

\pf (i) Since $W'$ is generated by $\iota_1^{'*}, \cdots,\iota_{l-1}^{'*}$, we may assume $w=\iota_j^{'*}$ for some $j \in \{1,\cdots,l-1\}$. A direct computation shows: \begin{align*}
    \phi(\iota_j^{'*}\cdot v) &= \iota_j^{'*}\cdot v-\frac{\inproc{\iota_j^{'*}\cdot v}{h_l}}{\inproc{h_l}{h_l}}h_l\\
     &= \iota^{*}_j\cdot v-\frac{\inproc{\iota^{*}_j\cdot v}{h_l}}{\inproc{h_l}{h_l}}h_l \quad\text{ since }(\iota^{*}_j-\iota_j^{'*})\cdot v \in \RR\cdot h_l\\
    &= \iota^{*}_j\cdot \left(v-\frac{\inproc{v}{h_l}}{\inproc{h_l}{h_l}}h_l\right) \quad\text{ by Lemma \ref{equi} and the fact $\iota^*_j\cdot h_l=h_l$.}\\
    &= \phi_W(\iota^{'*}_j)\cdot \phi(v).
\end{align*}

(ii) The inclusion of the image of $\Tbar'$ under $\phi$ in $\Tbar\cap h_l^\perp$ follows directly from (i). Conversely, consider any $v \in \Tbar\cap h_l^\perp$. Suppose $v=\lim_{m\to\infty} v_m$, where $\{v_m\}$ is a sequence of points in $\mathcal{T}$. By the continuity of $\phi$, we have $v=\phi(v)=\lim(\phi(v_m))$. Thus, it suffices to show that $\phi(v) \in \phi(\mathcal{T}')$ for any $v \in \mathcal{T}$. Let $w$ be a reduced word in $W$, and write $v=w\cdot(v'+u)$ with $v' \in \Tbar'$ and $u \in \RR_{\ge 0}h_l$. If $w$ contains no $\iota^*_l$, then 
$$\phi(v) = \phi\left(w\cdot v'+u\right) = \phi\left(w\cdot v'\right)\in \phi(\mathcal{T}').$$
Thus, we may assume $w=\ww\iota^*_lw'$, where $\ww \in W$ and $w' \in \langle\iota^*_j\rangle_{j=1,\cdots,l-1}$. Now, compute: 
\begin{align*}
    \phi(v) &= \phi\left(\ww\iota^*_lw'\cdot v'+\ww\iota^*_l\cdot u\right)\\
    &= \phi\left(\ww\iota^*_l(\phi_W^{-1}(w') \cdot v'+(w'-\phi_W^{-1}(w'))\cdot v')+\ww\iota^*_l\cdot u\right)\\
    &=\phi\left(\ww\phi_W^{-1}(w')\cdot v'\right)+\phi\left(\ww\iota^*_l\cdot((w'-\phi_W^{-1}(w'))\cdot v'+u)\right)
\end{align*}
By induction on the length of $w$, we may assume $\phi\left(\ww\phi_W^{-1}(w')\cdot v'\right) \in \phi(\mathcal{T'})$. Since $(w'-\phi_W^{-1}(w'))\cdot v' \in \RR_{\ge 0} h_l$ (Lemma \ref{lem}), it remains to show $\phi\left(\ww\iota^*_l \cdot h_l\right)\in \phi(\mathcal{T'})$. If $\ww=id$, the claim follows from Proposition \ref{iota}. For nontrivial $\ww$, using Proposition \ref{iota} again, we have, for $j\ne l$, as follows: \begin{align*}
    \iota^*_j\iota^*_l\cdot h_l &= \iota^*_j \cdot\left(-h_l+ \sum_{i\ne l} b_{il} h_i\right)\\
    &= -h_l+\sum_{i\ne j,l} b_{il} h_i +b_{jl}(\iota_j^{'*}\cdot h_j+b_{jl}h_l)\\
    &\in \iota_j^{'*}\left(\Abar'\right)+\RR_{\ge 0}h_l.
\end{align*}
By the inductive hypothesis, this implies that $\phi\left(\ww\iota^*_l \cdot h_l\right)\in \phi(\mathcal{T'})$, completing the proof of (ii).

(iii) Clearly, $\RR_{\ge 0}\cdot h_l \cup (\Tbar\cap h_l^\perp)$ is contained in $\Tbar\cap h_l^{\le 0}$. Conversely, for any $v\in \Tbar\cap h_l^{\le 0}$, we have $$v=\phi(v)+\frac{\inproc{v}{h_l}}{\inproc{h_l}{h_l}}h_l.$$ By the argument in (ii), we conclude that $\phi(v) \in \phi(\Tbar')$, completing the proof of (iii).\qed
\\
\begin{Lemma} \label{lem}
    Let $i\in\l$. We define a vector $v\in N^1(X)$ to be \emph{$i$-negative} if it satisfies $$\inproc{\alpha_i}{v}\le 0, \inproc{\alpha_j}{v} \ge 0, \text{ and } \inproc{\alpha_i}{v+\iota_j^*v}\ge 0$$ for all $j \in \l\setminus \{i\}$. Assume $v$ is $i$-negative, then $\iota_j^*v$ is $j$-negative for all $j\in \l\setminus \{i\}$. In particular, $wv$ is $i$-negative if $v\in \Abar$ and $w=\iota_i^*w_2\cdots w_m$ is a nontrivial reduced word in $W$.
\end{Lemma}
\pf Let $j \in \l\setminus \{i\}$. By hypothesis, we have $\inproc{\alpha_j}{\iota_j^*v}\le 0$. For $k\in \l\setminus \{i,j\}$, we compute: \begin{align*}
    &\inproc{\alpha_k}{\iota_j^*v} = \inproc{\alpha_k+b_{jk}\alpha_j}{v} \ge 0, \text{ and } \\
    &\inproc{\alpha_j}{\iota_j^*v+\iota_k^*\iota_j^*v}=\inproc{-\alpha_j}{v}+\inproc{b_{jk}\alpha_k+(b_{jk}^2-1) \alpha_j}{v}=\inproc{b_{jk}\alpha_k+(b_{jk}^2-2) \alpha_j}{v}\ge 0.
\end{align*} For the case when $k=i$, we observe that $$\inproc{\alpha_i}{\iota_j^*v}\ge-\inproc{\alpha_i}{v}\ge 0,$$ and \begin{align*}
    &\inproc{\alpha_j}{\iota_j^*v+\iota_i^*\iota_j^*v}=\inproc{b_{ij}\alpha_i+(b_{ij}^2-2) \alpha_j}{v} \ge 
    \inproc{(-\frac{b_{ij}^2}{2}+b_{ij}^2-2)\alpha_j}{v}\ge 0.
\end{align*} Therefore, we conclude that $\iota_j^*v$ is $j$-negative. It is straightforward to verify that $\iota_j^*v$ is $j$-negative for all $j\in\l$ and $v\in \Abar'$. Thus, the proof follows by induction on the length of $w$. \qed

\noindent\textbf{Decomposition of points in $\partial\Tbar$.}
The analog of \cite[Theorem 2.4.2]{FLT23} to the case $n > 1$ provides a decomposition of points on $\partial\Tbar$, which will play a crucial role in the proof of Theorem \ref{vol}.

\begin{Thm} \label{decomp} (cf. \cite[Theorem 2.4.2]{FLT23})
    Assume $nl \ge 4$, and $l \ge 3$ if $n=2$. For every $p\in \pTbar$, one of the following holds:
    \begin{enumerate}
        \item The direction of $p$ belongs to $\Lambda_{\l}$, namely, the image of $p$ under $\Tbar\longrightarrow\Tbar_+$ is in $\Lambda_{\l}$.
        \item There exist $i,j \in \l$ and $w\in W$, such that $w\cdot p =a_\lam v_\lam+\sum_{k\ne i,j}a_kh_k$ for some $a_\lam,a_k \ge 0$, where the direction of $v_\lam$ belongs to $\Lambda_{\{i,j\}}$ and $\inproc{v_\lam}{\alpha_k} > 0$ for $k\ne i,j$.
    \end{enumerate}
\end{Thm}
\pf This follows directly from Theorem \ref{boundary}. \qed

\begin{Rmk}
    See \cite[Theorem 2.4.2]{FLT23} for a stronger result in the case $n=1$.
\end{Rmk}

\subsection{Behavior of volume functions on $\partial\Tbar$}

\begin{Lemma} \label{compact}
    Suppose $n>1$, $nl \ge 4$, and $l \ge 3$ if $n=2$, then the set $\Fbar^1:=\{v\in \Fbar\mid\inproc{v}{v}=1\}$ is compact.
\end{Lemma}

\pf By (\ref{alpha}) and straightforward computation, we obtain
$$\omega_{\ihat j}:=\frac{b_{ij}}{2}\alpha_i+\alpha_j \in \sum_{k\ne i}\RR_{>0}h_k$$ for any $i\ne j \in \{1,\cdots,l\}$. In fact, we will show that all the extremal rays of $\Fbar$ are spanned by $\omega_{\ihat j}$. Using the definitions of $\Abar$ and $\Cbar$, it follows that $\Fbar$ is the intersection of the $2l$ hyperplanes
\begin{align*}
    &\{v \mid \inproc{v}{h_i}\ge0\} \text{ for }i=1,\cdots,l,\\
    &\{v \mid \inproc{v}{\alpha_j}\ge0\} \text{ for }j=1,\cdots,l.
\end{align*}
Observe that $$\inproc{\omega_{\ihat j}}{h_k}=\inproc{\frac{b_{ij}}{2}\alpha_i+\alpha_j}{h_k}\begin{cases}
    =0 \quad\text{ if }k \ne i,j,\\
    >0 \quad\text{ if }k = i,j,
\end{cases}$$ and $$\inproc{\omega_{\ihat j}}{\alpha_k}=\inproc{\frac{b_{ij}}{2}\alpha_i+\alpha_j}{\alpha_k}\begin{cases}
    =0 \quad\text{ if } k = i,\\
    >0 \quad\text{ if } k \ne i,
\end{cases}$$
Conversely, for any given $i,j$, if $v \in \Fbar$ satisfies $\inproc{v}{h_k}=\inproc{v}{\alpha_i}=0$ for $k\ne i,j$ must be a positive multiple of $\omega_{\ihat j}$. This proves that $\omega_{\ihat j}$ span extremal rays of $\Fbar$.

We then show that no additional extremal rays exist. Suppose $e \in \partial\Fbar$ spans an extremal ray. Then $e$ must belong to at least $l-1$ hyperplanes defined by $h_i^\perp$ or $\alpha_j^\perp$. Assume $e\in \displaystyle(\bigcap_{i \in S_h}h_i^\perp)\cap (\bigcap_{j \in S_h^c}h_{j}^{>0})$ for some subset $S_h \subset \l$ with $|S_h|\le l-1$. Consequently, $e$ lies in at least $l-1-|S_h|$ hyperplanes $\alpha_j^\perp$. Observe that $e \in h_j^\perp\cap\alpha_j^\perp$ implies $e=0$. Therefore, we know that $e \in \displaystyle\bigcap_{j\in S_h^c \setminus\{k\}}\alpha_{j}^\perp$ for some $k\in S_h^c$. Write $e=\sum_{j \in S_h^c} c_j\alpha_j$ for some positive coefficients $c_j$. Then, for $m \in S_h^c \setminus \{k\}$, we have $\sum_{j \in S_h^c} c_j \inproc{\alpha_j}{\alpha_m}=-\sum_{j \in S_h^c} c_jB_{jm}=0$. Assuming without loss of generality that $c_k=1$, it follows that $-B_{km}=\sum_{j \in S_h^c \setminus \{k\}} c_jB_{jm}$. Note that for $j\ne m$, the only possible values of $B_{jm}$ are $-n$ and $-n-\frac{1}{2}$ (see Remark \ref{bij remark}(b)). Therefore, we have the inequality $$n\le -B_{km}=\sum_{j \in S_h^c \setminus \{k\}} c_jB_{jm}\le c_m-\sum_{j \in S_h^c \setminus \{k,m\}} c_jn.$$ Adding these inequalities over $m \in S_h^c \setminus \{k\}$, we obtain $$n(|S_h^c|-1)\le\sum_{i\in S_h^c \setminus \{k\}}c_i-n(|S_h^c|-2)\sum_{i\in S_h^c \setminus \{k\}}c_i.$$ Thus, we conclude that $|S_h|\ge l-2$; otherwise, $\sum_{i\in S_h^c \setminus \{k\}}c_i\le 0$, which is a contradiction. If $|S_h|=l-1$, then $e$ is a positive multiple of $\alpha_k$ for some $k$. We find that $\inproc{e}{\alpha_k}<0$, which is again a contradiction. Therefore, $S_h^c = \{i,j\}$, and $\inproc{e}{\alpha_i}=0$ for some $i$ and $j$. This implies that $e$ is a positive multiple of $\om_{\ihat j}$.

We can now prove that $\Fbar^1$ is compact. For any $a \in \Fbar$, we can write $a$ as a linear combination of the extremal rays $\omega_{\ihat j}$ with nonegative coefficients. Let $a=\sum a_{ij} \omega_{\ihat j} \in \Fbar^1$, then $$\inproc{a}{a}=\sum_{i\ne j,k\ne l}a_{ij}a_{kl}\inproc{\omega_{\ihat j}}{\omega_{\khat l}} = 1.$$
Note that $$\inproc{\om_{\ihat j}}{\om_{\khat l}} = \frac{b_{ij}}{2}\cdot\frac{b_{kl}}{2}\cdot\frac{b_{ik}}{2}+\frac{b_{ij}}{2}\cdot\frac{b_{il}}{2}+\frac{b_{kl}}{2}\cdot\frac{b_{kj}}{2}+\frac{b_{jl}}{2}>0,$$ and hence, $$0<a_{ij}a_{kl} \le \left(\min{\{\inproc{\om_{\ihat j}}{\om_{\khat l}} \mid i \ne j, k\ne l, \text{ and } i,j,k,l \in \{1,\cdots,l\}\}}\right)^{-1}$$
In particular, $|a_{ij}|$ is bounded. This proves the compactness of $\Fbar^1$. \qed

\begin{Rmk}
    The extremal rays of $\Fbar$ in the case $n=1$ have been studied in \cite[Proposition 2.3.6]{FLT23}. However, in this case, for $n=1$, the set $\Fbar^1$ is not compact, and the asymptotic behavior of volume function near the extremal rays was explored in  \cite[\S 4.2]{FLT23}.
\end{Rmk}

\begin{Thm} \label{vol}
    Assume $n>1$, $nl \ge 4$, and $l \ge 3$ if $n=2$. Let $X$ be a general complete intersection Calabi-Yau subvariety in the product $\PP^{n} \times \cdots \times \PP^{n}$ of $l$ copies of $\PP^n$ defined as the intersection of $n$ ample divisors. Then, for any pseudoeffective divisor $p$ with $\vol_X(p)=0$, any ample $A$, and sufficiently small $s>0$, one of the following holds:
    \begin{enumerate}[label=(\alph*)]
        \item $\vol_X(p+sA)\asymp s^{\frac{kn}{2}}$ for some $k=1 \text{ or } l-1$,
        \item $\vol_X(p+sA)\asymp s^{kn}$ for some $k=1, \cdots, l-2$.
    \end{enumerate}
\end{Thm}

\pf (cf. \cite[Proposition 3.2.7 and \S 4.4]{FLT23}) Since the behavior of $\vol_X(p+sA)$ is independent of the choice of $A$ (\cite[Lemma 2.2(1)]{JW23}), we may assume $A$ belongs to the fundamental domain of $\CH$.
\\

\noindent {\it Case 1}: Assume that the direction of $p$ belongs to $\Lambda_{\l}$ ((i) of Theorem \ref{decomp}), then $\inproc{p+sA}{p+sA}=s^2\inproc{A}{A}+2s\inproc{p}{A}$, and $p+sA$ belongs to the $W$-orbit of $\Fbar$. Therefore, \begin{align*}
    \vol(p+sA)&=\inproc{p+sA}{p+sA}^{\frac{\dim X}{2}}\vol\left(\frac{p+sA}{\inproc{p+sA}{p+sA}^{\frac{1}{2}}}\right)\\
    &\asymp \inproc{p+sA}{p+sA}^{\frac{\dim X}{2}} \quad\text{ by Lemma \ref{compact} and the continuity of $\vol(\cdot)$}\\
    &\asymp s^{\frac{\dim X}{2}}\\
    &=s^{\frac{(l-1)n}{2}}.
\end{align*}

\noindent {\it Case 2}: Assume that $w\cdot p =a_\lam v_\lam+\sum_{k\ne i,j}a_kh_k$  as described in (ii) of Theorem \ref{decomp}. Without loss of generality, we may assume $w=1$. Let $A=\sum_{i=1}^la_i'h_i$ and $S=\{k\mid a_k\ne 0\}$. If $a_\lam =0$, then \begin{align*}
    \vol(p+sA)&=\vol\left(\sum_{s\in S} a_sh_s+s\sum_{i=1}^la_i'h_i\right)\\
    &\asymp s^{(\dim X-|S|)n}\\
    &=s^{(l-|S|-1)n},
\end{align*}
Notice that $1\le|S|\le l-2$. Therefore, $\vol(p+sA) \asymp s^{kn}$ for some $k=1,\cdots,l-2$.

Suppose now $a_\lam> 0$ and the direction of $v_\lam$ belongs to $\Lambda_{\{i,j\}}$. Let $a_\lam v_\lam=\sum_{k=1}^l x_kh_k$. Following the notation in \S\ref{preliminaries}, let $X'$ denote the restriction of $X$ obtained via the natural embedding of $\displaystyle\left(\prod_{k=i,j}\PP^n\right) \times \left(\prod_{k\in \l\setminus \{i,j\}}\{\text{point}\}\right)$ into $\displaystyle\prod_{k\in \l}\PP^n$. Define $A'=a_i'h_i+a_j'h_j$ and $p'=x_ih_i+x_jh_j$, then $p'$ is the a multiple of an element in $\Lambda'_{\{i.j\}} \subset \Tbar'$. Then we have\begin{align*}
    p+sA &= p'+sA'+\sum_{k\ne i,j}(x_k+a_k+sa_k')
\end{align*}
By Theorem \ref{decomp}, we know that $x_k>0$ for $k\ne i,j$. Hence, for sufficiently small $s$, $x_k+a_k+sa_k'$ is bounded above and below by positive constants when $k\ne i,j$. For the purpose of computing $\vol(p+sA)$, we select $w(s) \in W'$ such that $w(s)\cdot (p'+sA') \in \Abar'$, and hence, $\phi_W(w(s))\cdot (p+sA)\in \Abar$. This allows us to compute the volumes $$\vol_X(p+sA)=\vol_X(\phi_W(w(s))\cdot (p+sA)) \asymp \vol_{X'}(w(s)\cdot (p'+sA'))\asymp\vol_{X'}(p'+sA')$$ Combining this result with the analysis in \textit{Case 1}, we conclude that $\vol_X(p+sA)\asymp s^{\frac{\dim X'}{2}}= s^{\frac{n}{2}}$.

\begin{Rmk}
\begin{enumerate}
    \item The case when $n=1$ has been studied in \cite[Theorem 3.2.5]{FLT23}.
    \item The case where $X$ is a general complete intersection Calabi-Yau subvariety in $\PP^3\times\PP^3$ defined by ample divisors of degrees $(1,1)$, $(1,1)$ and $(2,2)$ has been discussed in \cite[\S 3]{Les19}.
\end{enumerate}
\end{Rmk}

\begin{Rmk}
    Note that the proof also implies the same result for abstract Lorentzian groups. More precisely, let $V$ be a vector space of dimension $l\ge 2$ with basis $\{h_i|i=1,\cdots,l\}$, and let $B$ be an $l\times l$ matrix as described in Proposition \ref{Lorentzian}. We define the involutions $\iota_i^*$ and the bilinear form using (\ref{hihj}) and (\ref{iota_alpha}). Furthermore, let $f$ be any homogeneous function of degree $(l-1)n$ on the cone generated by $\{h_i\}$. We can then extend $f$ by the action of $W$. Using the same argument, one can show that Theorem \ref{vol} holds for such $f$, assuming Theorem \ref{decomp}. In particular, when $n=2$, {\it Case 3} in the proof reduces to an abstract Lorentzian group with $l=2$, which remains valid by the same reasoning as in {\it Case 1}.
\end{Rmk}

\begin{Rmk}
    For a general complete intersection Calabi-Yau subvariety $X$ in $\PP^\nn:= \PP^{n_1} \times \cdots \times \PP^{n_l}$, as described in Theorem \ref{general vol theorem}, Proposition \ref{iota} implies that the asymptotic behavior of $\vol_X(p+sA)$ is determined by the components $p_j+sa_j$ for $n_j\in J$, where $p+sA=\displaystyle\sum_{i=1}^l (p_i+sa_i)\cdot h_i$. More precisely, let $X'$ be the restriction of $X$ under the natural embedding $\displaystyle\left(\prod_{i\in J}\PP^n\right) \times \left(\prod_{i\in J^c}\{\text{point}\}\right) \longrightarrow\displaystyle\prod_{i\in \l}\PP^n$, then $\vol_X(p+sA) \asymp \vol_{X'}(\displaystyle\sum_{j\in J}(p_j+sa_j)\cdot h_j')$. This proves Theorem \ref{general vol theorem}.
\end{Rmk}

\begin{Ex}
    Let $X$ be a general complete
intersection Calabi-Yau subvariety in $\PP^3\times\PP^3\times\PP^3$ defined by ample divisors of degrees $(1,1,2)$, $(1,2,1)$ and $(2,1,1)$. Consider the eigenvector $v= \frac{1}{14} (5 - 3 \sqrt{5})h_1+ \frac{1}{14} (5 + 3 \sqrt{5})h_2+ h_3$ associated with the eigenvalue $\frac{1}{2} (47 + 21 \sqrt{5})$ of $\iota_1^*\iota_2^*$. For sufficiently small $s>0$ and an ample divisor $A$, \begin{align*}
    &\vol_X(v+sA)\asymp s^3,\\
    &\vol_X(h_3+sA)\asymp s^3,\\
    &\vol_X(v+h_3+sA)\asymp s^{3/2}.
\end{align*}
\end{Ex}

\section{Examples of volumes for $\PP^2$-bundles of $X$} \label{Ch4}
Consider $\PP^{(1,1,2)}=\PP^1\times\PP^1\times\PP^2$ and let $X$ is a general divisor in $2H_1+2H_2+3H_3$. In this setting, we have $J=\{1,2\}$ and, by Proposition \ref{iota}, the corresponding involutions are given by
    $$\iota_1^* = \left[\begin{array}{ccc}
-1 & 0 & 0 
\\
 2 & 1 & 0 
\\
 3 & 0 & 1 
\end{array}\right]\text{,      }
\iota_2^* =\left[\begin{array}{ccc}
1 & 2 & 0 
\\
 0 & -1 & 0 
\\
 0 & 3 & 1 
\end{array}\right]$$
Recall that $\Nef(X)$ is the simplicial cone generated by $h_j$. Therefore, for $x,y,z\ge0$ , the volume function $\vol_X(x,y,z):= \vol_X(xh_1+yh_2+zh_3)$ is given by the top intersection \begin{align*}
    (xh_1+yh_2+zh_3)^3_X &= ((xH_1+yH_2+zH_3)^3\cdot(2H_1+2H_2+3H_3))_{\PP^{(1,1,2)}}\\
    &=18xyz+6xz^2+6yz^2.
\end{align*}

\begin{Lemma}\cite[Lemma 3]{BN16}
    Let $V$ be an irreducible projective variety of dimension $v$, $A_0,\cdots, A_r$ Cartier divisors on $V$. We set $$\mathcal{E} = \cO_V(A_0)\oplus \cdots \oplus \cO_V(A_r)$$ and consider the projective bundle $\PP(\mathcal{E})$ of dimension $d = v+r$. Then $$\vol_{\PP(\mathcal{E})}(\cO_{\PP(\mathcal{E})}(1)) = \frac{d!}{v!}\int_{\lambda_0+\cdots+\lambda_r,\,\lambda\ge0}\vol_V(\lambda_0 A_0+\cdots+\lambda_r A_r) d \lambda_1\cdots d \lambda_r.$$
\end{Lemma} 

Now, let $\textbf{w}=w_1w_2\cdots$ be a reduced word in $W=\langle\iota_1^*,\iota_2^*\rangle$. By a direct computation, we know that the cone $\textbf{w}(\Nef(X))$ is generated by \begin{align*}
    &\textbf{w}(h_1) = (-2 n +1)h_1+(2n)h_2+ (6 n^{2}-3 n)h_3,\\
    &\textbf{w}(h_2) = (-2 n)h_1+ (2 n +1)h_2+ (6 n^{2}+3 n)h_3,\\
    &\textbf{w}(h_3) = h_3 
\end{align*}if $|\textbf{w}| = 2n$ and $w_1=\iota_1^*$, and \begin{align*}
    &\textbf{w}(h_1) = (-2 n -1)h_1+(2 n +2)h_2+(6 n^{2}+9 n +3)h_3,\\
    &\textbf{w}(h_2) = (-2 n)h_1 +(2 n +1)h_2+(6 n^{2}+3 n)h_3\\
    &\textbf{w}(h_3) = h_3
\end{align*} if $|\textbf{w}| = 2n+1$ and $w_1=\iota_1^*$.

Let $S$ be the set $$\left(\bigcup_{\{\textbf{w} \mid w_1=\iota_1^*,|\textbf{w}|\ge 2\}}\textbf{w}(\Nef(X)) \right)\bigcap \{ax+by+z = k\}.$$ Suppose there exist integral points $A_0, A_1,A_2 \in \{ax+by+z = k\}$ such that their convex hull contains $S$ and does not intersect $\overline{\text{Eff}}(X) \setminus S$. Under this construction, the volume $\vol_X(\cO_X(1))$ is a rational multiple of \begin{align*}
    &V(a,b,k):=\int_S\vol(x,y,z) dS\\ &= \sum_{n=1}^\infty \left(\int_{(\iota_1^*\iota_2^*)^n\text{Nef}(X)\bigcap \{ax+by+z = k\}}\vol(x,y,z)dS+\int_{(\iota_1^*\iota_2^*)^n\iota_1^*\text{Nef}(X)\bigcap \{ax+by+z = k\}}\vol(x,y,z)dS\right)
\end{align*}

To simplify the computation, we choose $a,b$ such that the plane $ax+by+z = k$ intersects each $\textbf{w}(h_i)$ at exactly one point. A computation using Maplesoft then shows:

\begin{align*}
    &V(a,b,k)\\
    =& \frac{k^3}{10}\sum_{n=1}^\infty \frac{\left(3n^2 - (a - b)n + \frac{a}{4}+ \frac{b}{4} + \frac{3}{8}\right)}{\left(3n^2 - (a - b - \frac{3}{2})n + \frac{b}{2}\right)\left(3n^2 - (a - b + \frac{3}{2})n +\frac{a}{2}\right)}\\
    &+\frac{k^3}{10}\sum_{n=1}^\infty\frac{\left(3n^2 - (a - b - 3)n - \frac{a}{4} + \frac{3b}{4} + \frac{9}{8}\right)}{\left(3n^2 - (a - b - \frac{3}{2})n + \frac{b}{2}\right)\left(3n^2 - (a - b - \frac{9}{2})n - \frac{a}{2} + b + \frac{3}{2}\right)}\\
    =&\frac{k^3\left(4 a^2 + 4 b^2 - 8 a b - 12 a - 12 b - 9\right)}{120(\alpha-\beta)(a^2+b^2-2ab-3a-3b)}\left(\Psi\left(\frac{1}{2}-\beta\right)-\Psi\left(\frac{1}{2}-\alpha\right)+\Psi\left(1-\beta\right)-\Psi\left(1-\alpha\right)\right)\\
    &+\frac{k^3\left(2a^2+2b^2-4ab-3a-9b-18\right)}{20(a-2b-3)(a^2+b^2-2ab-3a-3b)}\\
\end{align*}
where $\alpha(a,b)>\beta(a,b)$ are the roots of $3x^2-(a-b-\frac{3}{2})x+\frac{b}{2}$.

In particular, $$V(-1,2,1)=\frac{1}{64}+\frac{\sqrt{33}}{264}\left(\Psi\left(\frac{7}{4}+\frac{\sqrt{33}}{12}\right)-\Psi\left(\frac{7}{4}-\frac{\sqrt{33}}{12}\right)+\Psi\left(\frac{5}{4}+\frac{\sqrt{33}}{12}\right)-\Psi\left(\frac{5}{4}-\frac{\sqrt{33}}{12}\right)\right).$$
We conclude Theorem \ref{Psi}.
\bibliographystyle{alphaurl} % or another style like alpha, apalike, IEEE, etc.
\bibliography{ref} % Use the name of your .bib file, without the .bib extension

\end{document}